\documentclass[12pt]{amsart}
\usepackage{amsfonts,amsmath,amsthm,oldgerm,amssymb,amscd}

\usepackage[all,cmtip]{xy}
\newcommand{\Q}{{\mathbb Q}}
\newcommand{\R}{{\mathbb R}}

\newcommand{\C}{{\mathbb C}}

\newcommand{\A}{{\mathbf A}}
\newcommand{\G}{{\mathbf G}}

\renewcommand{\H}{{\mathbb H}}

\newtheorem{theorem}{Theorem}

\theoremstyle{definition}

\numberwithin{definition}{section}

\numberwithin{equation}{section}

\theoremstyle{remark}
\newtheorem*{remark}{Remark}
\newtheorem*{acknowledgement}{Acknowledgement}

\begin{document}

\title[Conjugate non-homeomorphic varieties]{An example
  of non-homeomorphic  conjugate varieties}

\author{C.~S.~Rajan}

\address{Tata Institute of Fundamental  Research, Homi Bhabha Road,
Bombay - 400 005, INDIA.}  \email{rajan@math.tifr.res.in}


\begin{abstract}
We give examples of smooth quasi-projective varieties over complex
numbers, in the context of connected Shimura varieties, 
which are not homeomorphic to a conjugate of itself by an automorphism
of the complex numbers.

\end{abstract}
 
\maketitle

\section{Introduction}  
Let $X$ be a quasi-projective variety defined over $\C$. Suppose 
$\sigma$ is  an automorphism of $\C$.  Denote by 
$X^{\sigma}:=X\times_{\sigma}\C$, the conjugate of $X$ by the
automorphism $\sigma$ of $\C$, obtained by applying the automorphism
$\sigma$ to the coefficients of the polynomials defining $X$. 
It is known that the varieties $X$ and $X^{\sigma}$ have the same
Betti numbers. In \cite{Se}, Serre gave an example such that the
topological spaces $X(\C)$ and $X^{\sigma}(\C)$ are not
homeomorphic.

 Recently, Milne and Suh \cite{MS}, gave further examples
in the context of connected Shimura varieties. Their method is to
find a conjugate such that the reductive group underlying the Shimura
datum is different, and then apply the super-rigidity results of
Margulis. 

Our examples are in the same context as  that of Milne and Suh, but we
work with Shimura's construction of canonical models 
(\cite{Sh}). Shimura's construction allows us to identify the
adelic congruence subgroup defining the conjugate variety as a
conjugate by an element of the adjoint group. We then appeal to Mostow
rigidity and the failure of strong approximation (or
non-triviality of class number) for the adjoint group
to get at the desired examples. In our example, the congruent lattices
defining the variety and it's conjugate  are commensurable.  Earlier
in \cite{R}, we observed using Shimura's construction 
coupled with the theorems of Labesse and
Langlands on the mulitplicity of cusp forms for $SL(1, D)$, 
 that a Galois twist of these spaces attached to $SL(1, D)$ 
over the reflex field
preserves the spectrum of the Laplacian; this provides    
examples of locally symmetric spaces
attached to a quaternion division algebra over a number field
which are isospectral but
not isometric,

\section{The example}
Let $F$ be a totally real number field of degree at least two. Let $D$
be an indefinite  quaternion  algebra defined over $F$. 
We assume that $D$ is split
at exactly one real place, say $\tau_1$  of $F$. This assumption
allows us to assume that the reflex field of $(F, \tau_1)$  to be $F$
itself.   Let $V$ be a vector
space of rank $n\geq 2$ over $D$, equipped with a hermitian inner product
with respect to the standard involution on $D$. We assume that the
inner product is definite on the spaces $V\otimes_{\tau}\R$, for all
real embeddings $\tau$ of $F$ different from $\tau_1$. In particular, 
since we have assumed that the degree of $F$ is at least two, the form
$h$ is anisotropic. Let $G$ be the
group of unitary similitudes of $h$. We consider $G$ as an algebraic
group defined over $\Q$, and let $G_d$ the derived group of
$G$. We let $PG$ denote the projective group attached to $G$, the group
obtained by taking the quotient of $G$ modulo it's centre. 
Under our assumptions, it follows that 
\[ G_d(\R) \simeq Sp(2n, \R)\times ~\mbox{a compact group}, \quad
n\geq 2.\]

Let $K_{\infty}$ be a maximal compact subgroup of
$G_d(\R)$ and let $X=G_d(\R)/K_{\infty}$ be the
non-compact symmetric space associated to $G$. By our assumptions, 
 $ X$ is isomorphic to 
 the Siegel  upper half space  $\H_n$ of dimension $n$. 
Denote by  $\A$ the
adele ring of $F$, and by $\A_f$ the subring of finite adeles. 
Let $K$ be a
compact open subgroup of $G(\A_f)$, and let $K_d=K\cap G_d(\A_f)$. 
 Denote by    
\[\Gamma_K=G(\R)K\cap G(\Q) \quad \mbox{and}\quad
\Gamma_{d,K}=G _{d}(\R)K_d\cap G_d(\Q),  \]
the corresponding  arithmetic lattices in $G(\R)$ and
$G_d(\R)$ respectively.  
We assume that $K$ is such that
$\Gamma_{d,K}$ is  torsion-free, and the
natural inclusion $\Gamma_{d,K}\subset \Gamma_K $ is an
isomorphism modulo the centre of $\Gamma_K$.  

By a theorem of Baily-Borel, 
the quotient space
$X_K=\Gamma_K\backslash X$ is a connected,  smooth, projective
variety. The fundamental group  $\overline{\Gamma}_K$ of the variety
$X_K$ can be identified  with  the
projection of ${\Gamma_K}$ to
$PG(\R)$, and also with the lattice $\Gamma_{d,K}$ contained in
$G_d(\R)$.

For an element $x\in G(\A_f)$, denote by $K^x$ the conjugate
lattice $x^{-1}Kx$, and by $\overline{x}$ its image in $PG(\A_f)$. 
Further, let
  $N(\overline{K})$ denote the normalizer of $\overline{K}$ in
  $PG(\A_f)$, where $\overline{K}$ is the image of $K_d$ in
  $PG(\A_f)$.  The desired example is provided by the following
theorem: 
\begin{theorem} \label{adelicconjtheorem}
With notation and assumptions as above,    
suppose $x$ is an element in $G(\A_f)$ such
  that $\overline{x}$ does not belong to the set $N(\overline{K})PG(\Q)$.
  Then $X_K$ and $X_{K^x}$ are conjugate by an automorphism $\sigma$
  of $\C$, but the respective fundamental groups $\overline{\Gamma}_K$ and
  $\overline{\Gamma}_{K^x}$ are not isomorphic.  
 In particular,  $X_K$ and $X_{K^x}$
  are not homeomorphic.
\end{theorem}

\begin{remark} The hypothesis can be seen to hold from two different
  but related aspects of the arithmetic of algebraic groups. The 
normalizer $N(\overline{K})$ is a compact open subgroup of $PG(\A_f)$.
By the failure of strong approximation for the adjoint group
$PG$ (see \cite[Proposition 7.13]{PR}), the rational points $PG(\Q)$ 
are not dense in $PG(\A_f)$. 
Hence,  the hypothesis that $\overline{x}$ does not belong to the
double coset $N(\overline{K})PG(\Q)$ is satisfied 
provided $N(\overline{K})$ is
small enough. 

On the other hand, the adjoint group $PG$ is not simply
connected, hence has a non-trivial fundamental group. 
 The results of Section 8.2 
of \cite{PR}, show that the class group  of $G$ is non-trivial for
suitably chosen congruence lattices  in $PG(\A_f)$. This allows us to
work with large congruence lattices in $PG(\A_f)$. 
\end{remark}

\begin{proof}
We first show that the varieties $X_K$ and $X_{K^x}$ are conjugate by
an automorphism of $\C$. For this, 
 we recall Shimura's theory of canonical
models \cite{Sh}. 
Let $\nu: {G}\to \G_m$ be  the reduced norm. 
By class field theory,  the  subgroup $F^*\nu(K)$ of the idele
group $\A^*$   
defines an abelian
extension $F_K$ of $F$.  The reciprocity morphism of class field,
\[ \mbox{rec}: \A^*/F^*\to \mbox{Gal}(F^{ab}/F),\]
defines an element $\sigma(x)\in  \mbox{Gal}(F^{ab}/F)$ by the
prescription 
\[ \sigma(x)= \mbox{rec}(\nu(x)^{-1}).\]
As a consequence of the main theorem of
canonical models in \cite[Theorem 2.5, page 159, Section 2.6]{Sh}, the
variety $X_K$ has a model defined over the field $F_K$, and 
\begin{equation}\label{conjcanon}
 X_K^{\sigma(x)}\simeq X_{K^x}.
\end{equation}
Thus the varieties $X_K$ and $X_{K^x}$ are conjugate.

Suppose on the contrary, that $X_K$ and $X_{K^x}$ have isomorphic
fundamental groups. Since these spaces are Eilenberg-Maclane spaces,
there exists a homotopy equivalence
\[\phi: X_K\to X_{K^x}.\]
Since the lattices are irreducible in $PG(\R)$ and the real rank of
$PG$ is at least two, 
by Mostow rigidity \cite{Mo}, 
the spaces $X_K$ and $X_{K^x}$ are isometric. 

Hence, there exists $\overline{g}\in PG(\R)$ such that
\[ \overline{g}^{-1}\overline{\Gamma}_{K^x}\overline{g}=\overline{\Gamma}_K.\]
 Since the lattices $\overline{\Gamma}_K$ and
$\overline{\Gamma}_{K^x}$ are arithmetic and commensurable, it follows  by a
theorem of Borel (\cite{Bo}), that $\overline{g}\in PG(\Q)$. Hence
 there is an element $g\in G(\Q)$ satisfying,
\[  g^{-1}\Gamma_{d,K^x}g=\Gamma_{d,K}.\]
Consider now $G_d(\Q)$ embedded diagonally in $G_d(\A_f)$. By the strong approximation theorem
for $G_d$, the closure of $\Gamma_{d,K}$ in $G_d(\A_f)$ can be identified
with $K_d$. Further, the closure of $\Gamma_{d,K^x}$ in $G_d(\A_f)$
can be identified with  $g^{-1}K_d^xg$, where we now consider $g\in
G(\Q)$ as diagonally embedded in $G(\A_f)$. Hence,  we have 
\[ g^{-1}K_d^xg=K_d.\]
Projecting to $PG$, we obtain
\[
\overline{g}^{-1}\overline{x}^{-1}\overline{K}\overline{x}\overline{g}
=\overline{K}, \]
where $\overline{K}$ denotes the image of $K_d$ in $PG(\A_f)$.  This
implies that $\overline{x}\in N(\overline{K})PG(\Q)$, contradicting our
choice of $\overline{x}$.
\end{proof}

\begin{acknowledgement}
My sincere thanks to Dipendra Prasad for useful suggestions and
discussions. 
\end{acknowledgement}

\end{document}